\documentclass[12pt]{elsarticle}

\usepackage{latexsym,amsfonts,amsmath,theorem,amssymb}
\usepackage{graphicx}

\def\be{\begin{equation}}
\def\bea{\begin{eqnarray*}}
\def\ee{\end{equation}}
\def\eea{\end{eqnarray*}}
\def\ba{\begin{array}}
\def\ea{\end{array}}
\def\bi{\begin{itemize}}
\def\ei{\end{itemize}}
\newtheorem{theo}{Theorem}

\def\ZZ{{{\rm Z}\kern-.4em{\rm Z}}}
\def\RR{{{\rm I}\kern-.2em{\rm R}}}
\def\NN{{{\rm I}\kern-.2em{\rm N}}}
\def\TT{{{\rm T}\kern-.5em{\rm T}}}
\def\CC{{{\rm I}\kern-.5em{\rm C}}}

\journal{XXX}

\begin{document}
\begin{frontmatter}

\title{Estimates for Generalized Sparse Grid Hierarchical Basis Preconditioners}
\author[bonn]{Peter Oswald}
\ead{agp.oswald@gmail.com}

\address[bonn]{Institute for Numerical Simulation (INS),
University of Bonn,
Wegelerstr. 6-8,
D-53211 Bonn
}
\date{}
\begin{abstract}
We reconsider some estimates from the 1994 paper \cite{GO1994} concerning the hierarchical basis preconditioner
for sparse grid discretizations.
The improvement is in three directions: We consider arbitrary space dimensions $d>1$, give bounds  for generalized sparse grid spaces with arbitrary monotone index set $\Lambda$, and show that the bounds are sharp up to constants depending 
only on $d$, at least for a subclass of $\Lambda$ containing full grid, standard sparse grid spaces, and energy-norm optimized sparse grid spaces.
\end{abstract}
\end{frontmatter}

\section{Notation}\label{sec1} We consider the hierarchical basis (HB)  preconditioner for generalized sparse grid discretizations for generic $H_0^1(I^d)$-elliptic problems ($I^d$ is the $d$-dimensional unit cube) which has been analyzed for $d=2,3$ in \cite{GO1994}. The underlying hierarchical basis is a finite collection of dyadic blocks of the tensor-product Faber-Schauder system on $I^d$ (detailed definitions will be given in the next section). The HB preconditioner considered in this note  is different from the HB preconditioners of Bank, Dupont, and Yserentant (see \cite{B1997,Y1993}) which correspond to  isotropic multivariate versions of the Faber-Schauder system. The HB preconditioner for sparse grid discretizations is not optimal, nor suboptimal (better preconditioners have already been discussed  in \cite{GO1994,O1994} and more recently in \cite{GHO2014}). Nevertheless, due to the popularity of sparse grid methods we find it worthwhile to have a complete understanding of its properties.

In this note, we consider generalized sparse grid spaces $S_\Lambda$ for arbitrary dimension $d>1$ generated by
a monotone index set $\Lambda\subset \mathbb{N}^d$  which comes with a direct splitting
$$
S_\Lambda = \sum_{\beta\in\Lambda} S_\beta
$$
where each $S_\beta$ is the span of nodal basis functions associated with the dyadic block of the tensor-product Faber-Schauder system with index $\beta\in\mathbb{N}^d$. Correspondingly, each $v_{\Lambda}\in S_\Lambda$ has a unique decomposition
$$
v_\Lambda =\sum_{\beta\in\Lambda} s_\beta,\qquad s_\beta\in S_\beta.
$$
The associated Galerkin discretization for a generic $H_0^1(I^d)$-elliptic problem using the finite section of the
tensor-product  Faber-Schauder system associated with $S_\Lambda$ leads, after appropriate diagonal scaling, to a symmetric positive-definite algebraic system 
$$
A_{\Lambda} x_\Lambda = b_{\Lambda},
$$
with stiffness matrix $A_\Lambda$, right-hand side $b_\Lambda$, and solution vector $x_\Lambda$ representing the HB coefficients of the Galerkin projection of the $H^1_0$-elliptic problem onto $S_\Lambda$. We seek as good as possible estimates of the spectral bounds $\lambda_{\max}(A_\Lambda)$ and $\lambda_{\min}(A_\Lambda)$, and, consequently, of the spectral condition number
$$
\kappa_{S_\Lambda,HB}:=\kappa(A_\Lambda)=\frac{\lambda_{\max}(A_\Lambda)}{\lambda_{\min}(A_\Lambda)}.
$$ 
This is equivalent to estimating the stability bounds of the direct space splitting
\be\label{SS}
\{S_\Lambda;(\cdot,\cdot)_{H^1_0}\} =\sum_{\beta\in \Lambda} 
\{S_\beta;2^{2|\beta|_\infty}(\cdot,\cdot)_{L_2}\}.
\ee
This standard fact from additive Schwarz theory is explained in \cite{GO1994}. In particular, estimating $\lambda_{\max}(A_\Lambda)$ (up to generic constants depending on the ellipticity constants and $d$) is equivalent to finding the best constant $C_\Lambda$ in the inequality
\be\label{UE}
\|\sum_{\beta\in \Lambda} s_\beta\|^2_{H_0^1} \le C_\Lambda  \sum_{\beta\in \Lambda}
2^{2|\beta|_\infty}\|s_\beta\|^2_{L_2},
\ee
valid for all $s_\beta\in S_\beta$ with $\beta \in\Lambda$. Upper and lower estimates for $C_\Lambda$ are obtained in Section \ref{sec2}, they are  matching up to constants depending  on $d$ but not on $\Lambda$.

Estimates for $\lambda_{\min}(A_\Lambda)$ require bounds for the best constant $c_\Lambda$ in
the inequality 
\be\label{LE}
c_\Lambda \sum_{\beta\in \Lambda}
2^{2|\beta|_\infty}\|s_\beta\|^2_{L_2}\le \|\sum_{\beta\in \Lambda} s_\beta\|^2_{H_0^1}, \qquad s_\beta\in S_\beta,
\ee
opposite to (\ref{UE}), and will be given in Section \ref{sec3}. In the final section we summarize the results and show that the resulting condition number estimates for HB preconditioners are asymptotically sharp for certain families of $S_\Lambda$ with $d>1$ arbitrarily fixed, including the full grid spaces $V_k$
and standard sparse grid spaces $S_k$ if $k\to \infty$ .

\section{Notation and auxiliary facts}\label{sec11}
\subsection{Faber-Schauder functions} 
Denote by $T_k$ the uniform dyadic partition of $I=[0,1]$ of step-size $2^{-k}$, $k\in \mathbb{N}$. 
The univariate Faber-Schauder system on the unit interval (obeying zero boundary conditions) 
consists of dyadic shifts and dilates of the unit hat function
$\phi(t)=(1-|t|)_+$, $t\in \mathbb{R}$, defined blockwise as follows: For $k=1,2,...$, define the $k$-th dyadic block of the Faber-Schauder system as the collection of $2^{k-1}$ hat functions
$$
\phi_{2^{k-1}+i}(t):=\phi(2^kt-(2i+1)),\quad t\in [0,1],\qquad i=0,1,\ldots, 2^{k-1}-1,
$$
with non-overlapping support. These are the standard nodal basis functions of the linear finite element space over $T_k$ associated with the interior nodal points of $T_k$ of the form $(2i+1)2^{-k}$. The index set for this block is denoted $J_k$.
The union of all blocks is the univariate Faber-Schauder system (with zero boundary conditions).

The multivariate tensor-product Faber-Schauder system is defined as the collection of  the functions
$$
\phi_\alpha(x):=\prod_{i=1}^d \phi_{\alpha_i}(x_i),\qquad \alpha=(\alpha_1,\ldots,\alpha_d)\in \mathbb{N}^d.
$$
We organize it into blocks associated with the multi-index sets $J_\beta:=J_{\beta_1}\times\ldots \times J_{\beta_d}$,
where $\beta\in\mathbb{N}^d$.
The Faber-Schauder functions $\phi_\alpha$ in the same block have non-overlapping support, and are shifts of each other.
Together with their tensor-product structure this allows us to obtain explicit formulas for their $L_2$ and $H^1_0$ norms
(for convenience, the latter is defined as 
$$
\|u\|_{H_0^1}^2 := \int_{I^d} |\nabla u|^2\, dx,\qquad u\in H^1_0).
$$
More precisely,
\be\label{phiN}
\|\phi_\alpha\|_{L_2}^2 = \frac{2^d}{3^{d}}2^{-|\beta|_1},\qquad 
\|\phi_\alpha\|_{H^1_0}^2 = \frac{2^d}{3^{d-1}} 2^{-|\beta|_1}\sum_{i=1}^d 2^{2\beta_i},\qquad \alpha\in J_\beta.
\ee
We use the notation $|\beta|_1:=\sum_{i=1}^d \beta_i$ and  $|\beta|_\infty=\max_{i=1,\ldots,d} \beta_i$.

\subsection{Discretization spaces and decomposition norms}
Let $S_\beta$ denote the finite-dimensional space spanned by all $\phi_{\alpha}$ with $\alpha\in J_\beta$. Then
$$
V_\beta = \sum_{\beta'\le \beta} S_{\beta'}, \qquad \beta\in \mathbb{N}^d,
$$
is a direct sum splitting of the anisotropic full-grid space $V_\beta$ of all $d$-linear finite element functions over the anisotropic tensor-product partition $T_\beta:=T_{\beta_1}\times \ldots\times T_{\beta_d}$ (by $\beta'\le \beta$ we mean that $\beta'_i\le \beta_i$ for all $i=1,\ldots,d$). The full grid space $V_k$ refers to the isotropic case $\beta=(k,\ldots,k)$ needed in approximation schemes related to uniform grid refinement. 

In addition to the subspace families $\{S_\beta\}$ and $\{V_\beta\}$,
we also need the family
$$
W_\beta := V_\beta \ominus (\sum_{\beta'<\beta} V_{\beta'}) = \otimes_{i=1}^d (V_{\beta_i}\ominus V_{\beta_i-1})
$$
of $L_2$ orthogonal subspaces which provides splittings of the various spaces of interest to us (the convention is $V_0=\{0\}$). In particular,
\be\label{W}
L_2(I^d)= \oplus_{\beta\in\mathbb{N}^d} W_\beta,\qquad V_\beta = \oplus_{\beta'\le \beta} W_{\beta'}.
\ee

The main focus is the investigation of the generalized sparse grid space
\be\label{S}
S_\Lambda := \mathrm{span}(\{V_\beta:\,\beta\in \Lambda\}) =  \mathrm{span}(\{W_\beta:\,\beta\in \Lambda\})= \mathrm{span}(\{S_\beta:\,\beta\in \Lambda\}),
\ee
where $\Lambda\subset \mathbb{N}^d$ is a monotone set, i.e., $\beta\in \Lambda$ implies $\beta'\in\Lambda$ for all $\beta'\le \beta$. The standard sparse grid spaces $S_k$ correspond to the choice $\Lambda=\{\beta\in \mathbb{N}^d:\; |\beta|_1\le k+d-1\}$.
Note that the anisotropic full grid spaces $V_\beta$ are also special instances of the family of $S_\Lambda$ spaces (take $\Lambda=\{\beta'\in\mathbb{N}^d:\;\beta'\le \beta\}$). 

Each $v_\Lambda\in S_\Lambda$ can be non-uniquely decomposed as
$$
v_\Lambda = \sum_{\beta\in\Lambda} v_\beta,\qquad v_\beta\in V_\beta,
$$
but also uniquely represented with respect to either the $S_\beta$ or $W_\beta$ subspaces of $V_\Lambda$:
\be\label{vLa}
v_\Lambda = \sum_{\beta\in\Lambda} s_\beta = \sum_{\beta\in\Lambda} w_\beta,\qquad s_\beta\in S_\beta,\quad w_\beta\in W_\beta.
\ee
Each of these representations of $v_\Lambda\in S_\Lambda$ has its own merits. In particular, $L_2$ orthogonal representations (also called prewavelet (PW) representations) with respect to the 
$W_\beta$ subspaces allow us to effectively express $L_2$ and $H_0^1$ norms. For any $v_\Lambda\in S_\Lambda$ we have the identity
$$
\|v_\Lambda\|_{L_2}^2 = \sum_{\beta\in \Lambda} \|w_\beta\|_{L_2}^2,
$$
and the two-sided norm equivalence
\be\label{PWN}
c_{PW}\|v_\Lambda\|_{PW}^2\le \|v_\Lambda\|_{H^1_0}^2
\le C_{PW}\|v_\Lambda\|_{PW}^2,\qquad \|v_\Lambda\|_{PW}^2 := \sum_{\beta\in \Lambda} 2^{2|\beta|_\infty}\|w_\beta\|_{L_2}^2,
\ee
which holds with constants $0<c_{PW}\le C_{PW}<\infty$ depending on $d$, only. We also need an inequality for arbitrary decompositions with respect to the isotropic full grid spaces $V_k$, namely 
\be\label{BPX}
\|\sum_{k=1}^K v_k\|_{H_0^1}^2 \le C_{BPX} \sum_{k=1}^K 2^{2k}\|v_k\|^2_{L_2},\qquad v_k\in V_k,
\ee
which holds with a constant depending on $d$. It is related to the BPX preconditioner, and can be obtained from (\ref{PWN}). We refer to \cite{O1994,Y1993} for more details on the inequalities (\ref{PWN}) and (\ref{BPX}).

In the remainder of this paper, we will be interested in 
identifying the constants in a similar two-sided estimate associated with the decomposition with respect to the $S_\beta$ spaces,
namely for comparing the $H_0^1$ norm with the HB norm
\be\label{HBN}
\|v_\Lambda\|_{HB}^2 := \sum_{\beta\in \Lambda} 2^{2|\beta|_\infty}\|s_\beta\|_{L_2}^2,
\ee
instead of the PW norm in (\ref{PWN}). As will be demonstrated in the next two sections, these constants also depend on characteristics of $\Lambda$, and not only on $d$. 

\subsection{Norms of some FE functions} 
We start with stating an immediate consequence of (\ref{phiN}) and the non-overlapping support property
of the nodal basis functions spanning $S_\beta$:
\be\label{SN}
\|s_\beta\|_{L_2}^2 =  \frac{2^d}{3^{d}} 2^{-|\beta|_1} \sum_{\alpha\in J_\beta} c_\alpha^2,\qquad s_\beta =\sum_{\alpha\in J_\beta} c_\alpha \phi_\alpha \in S_\beta.
\ee
A similar equality holds for the $H^1_0$ norm but we do not need it.

Next, we estimate the HB norm of the tensor-product hat function
$$
\psi_\beta(x) = \prod_{i=1}^d \phi(2^{\beta_i}x-2^{\beta_i-1}),\qquad  \beta\in \mathbb{N}^d.
$$
This is the nodal basis function in $V_\beta$ associated with the center of the cube $I^d$ (which is not in the tensor-product Faber-Schauder system 
unless $\beta=(1,\ldots,1)$). Since it is the tensor product of $d$ univariate hat functions $\phi(2^{\beta_i}x-2^{\beta_i-1})$ associated with the nodal point $1/2$, its HB decomposition is the tensor product of the univariate HB decompositions of the latter. It is easy to see that
the univariate HB decompositions are implied by the formula
$$
\phi(2^{m}t-2^{m-1})=\phi_{2^{0}}(t) - \sum_{l=2}^m \frac12(\phi_{2^{l-1}+2^{l-2}-1}(t) +\phi_{2^{l-1}+2^{l-2}}(t)),\quad t\in [0,1],
\quad m\in \mathbb{N},
$$
by setting $m=\beta_i$ and $t=x_i$, $i=1,\ldots,d$. Thus, if the multi-index $\beta'\le \beta$ has $r$ components $\beta'_i>1$ and
$d-r$ components $\beta'_i=1$ then the HB block $s_{\beta'}$ of $\psi_\beta$ is the sum of $2^r$ nodal basis functions with coefficient
$(-1/2)^r$. Thus, by (\ref{SN})
$$
\|s_{\beta'}\|_{L_2}^2 =  \frac{2^d}{3^{d}}2^{-|\beta'|_1} 2^r (1/2)^{2r} \ge 3^{-d} 2^{-|\beta'|_1}, \qquad \beta'\le \beta.
$$
Here, and throughout the paper, we denote by $c,C>0$ generic positive constants depending only on $d$ (which may be different in different places). Moreover, we use the notation $A\approx B$ if $c\,A \le B\le C\,A$. Substitution into the expression for the HB norm gives
\be\label{phiHB}
\|\psi_\beta\|_{HB}^2 = \sum_{\beta'\le\beta} 2^{2|\beta'|_\infty}\|s_{\beta'}\|_{L_2}^2\ge c\sum_{\beta'\le\beta} 2^{2|\beta'|_\infty -|\beta'|_1} \ge  c2^{|\beta|_\infty}.
\ee
The last inequality follows since among the indices $\beta'\le \beta$ there is at least one with the property 
$|\beta'|_1-(d-1)=|\beta'|_\infty=|\beta|_\infty$. There is also a matching upper bound 
$\|\psi_\beta\|_{HB}^2  \le  C2^{|\beta|_\infty}$ but we do not need it in the sequel.

Finally, we consider a different construction (a kind of lacunary HB series representation) which we need in the next section. For the following estimates to hold the set $\Lambda$ can be arbitrary, i.e., not necessarily monotone. Let
$$
\bar{s}_\beta = \sum_{\alpha\in J_\beta} \phi_\alpha\in S_\beta, \qquad \beta\in \mathbb{N}^d,
$$
and set
$$
\bar{s}_\Lambda = \sum_{\beta\in \Lambda} \bar{s}_\beta\in S_\Lambda.
$$
The functions $\bar{s}_\beta$ are tensor products of their univariate counterparts 
$$
\bar{s}_k(t) = 2^{k} \int_0^t r_k(\xi)\,d\xi, \qquad t\in [0,1],\qquad k=1,2,\ldots,
$$
where $r_k(t)=\mathrm{sign}(\sin(2^k\pi t))$ denotes the univariate Rademacher functions.
Below we need that the shifted functions
$$
\tilde{s}_k(t)= \bar{s}_k(t)-1/2, \qquad k=1,2,\ldots,
$$
form an orthogonal system in $L_2([0,1])$ (and are orthogonal to constants as well), with $L_2$ norm given by
$$
\|\tilde{s}_k\|^2_{L_2}=\|\bar{s}_k\|^2_{L_2}-\|\bar{s}_k\|^2_{L_1}+\frac14 = \frac13-\frac12+\frac14 =\frac1{12}.
$$
Here we have used the case $d=1$ of  the identity
$$
\|\bar{s}_\beta\|_{L_2}^2 = \frac{2^d}{3^d} 2^{-|\beta|_1} \sum_{\alpha\in J_\beta} 1 = \frac{2^d}{3^d} 2^{-|\beta|_1} 2^{|\beta|_1-d} =3^{-d},
$$
which is a consequence of (\ref{SN}).
From the last equality, we compute the HB norm of $s_\Lambda$ as
\be\label{Lac1}
\|\bar{s}_\Lambda\|_{HB}^2 = \sum_{\beta\in \Lambda} 2^{2|\beta|_\infty}\|\bar{s}_\beta\|_{L_2}^2=3^{-d} \sum_{\beta\in\Lambda} 2^{2|\beta|_\infty}. 
\ee

For the $L_2$ norm of $\bar{s}_\Lambda$ we can obtain the lower bound
\be\label{Lac2}
\|\bar{s}_\Lambda\|_{L_2}^2 \ge 4^{-d}|\Lambda|^2.
\ee
Indeed, since 
$$
\bar{s}_\beta(x)=\prod_{i=1}^d (\tilde{s}_{\beta_i}(x_i)+\frac12),
$$
using the orthogonality properties of the system $\{\tilde{s}_k(t)\}$ mentioned before, namely that
$$
\int_0^1 (\bar{s}_{k}(t)+\frac12)(\bar{s}_{k'}(t)+\frac12)\, dt = \frac{\delta_{kk'}}{12}+\frac14 \ge \frac14,
$$
we have
\bea
\|\bar{s}_\Lambda\|_{L_2}^2 &\ge& \sum_{\beta\in \Lambda} \sum_{\beta'\in \Lambda} 
  \prod_{i=1}^d \int_0^1 (\bar{s}_{\beta_i}(x_i)+\frac12)(\bar{s}_{\beta'_i}(x_i)+\frac12)\, dx_i\\
	&\ge& \sum_{\beta\in \Lambda} \sum_{\beta'\in \Lambda} 4^{-d} = 4^{-d}|\Lambda|^2.
\eea

\section{Upper estimates: $C_{\Lambda}$ }\label{sec2} 
We follow the approach adopted  for $d=2$ in \cite{GO1994}. Take an arbitrary $v_\Lambda\in S_\Lambda$, and consider its unique decomposition (\ref{vLa}) into HB blocks $s_\beta\in S_\beta$.
Define a partition of $\Lambda$ into non-overlapping index subsets 
\be\label{Lk}
\Lambda_k:= \{\beta \in \Lambda: \;|\beta|_\infty =k\},\qquad k = 1,\ldots,k_\Lambda,\quad k_\Lambda:=\max_{\beta\in\Lambda} |\beta|_\infty,
\ee
and gather the $s_\beta$ into blocks associated with these $\Lambda_k$. Obviously,
$$
v_k:= \sum_{\beta\in \Lambda_k} s_\beta \in V_k, \qquad k=1,\ldots,k_\Lambda,
$$
and according to (\ref{BPX})
$$
\|v_\Lambda\|^2_{H_0^1}=\|\sum_{k=1}^{k_\Lambda} v_k\|^2_{H_0^1} \le C_{BPX} \sum_{k=1}^\infty
2^{2k}\|v_k\|^2_{L_2}.
$$
If $|\Lambda_k|$ denotes the number of indices in $\Lambda_k$ then by the Cauchy-Schwarz inequality
$$
\|v_k\|^2_{L_2}=\|\sum_{\beta\in \Lambda_k} s_\beta\|^2_{L_2}\le |\Lambda_k|
\sum_{\beta\in \Lambda_k} \|s_\beta\|^2_{L_2},
$$
and, since $k=|\beta|_\infty$ for all $\beta\in\Lambda_k$, after substitution we see that the best constant in (\ref{UE}) satisfies
\be\label{CLambda1}
C_\Lambda \le C_{BPX} \cdot n_\Lambda,\qquad n_\Lambda:=\max_{1\le k\le k_\Lambda} |\Lambda_k|.
\ee

A matching lower bound for the best possible $C_\Lambda$ in (\ref{UE}) is given by the following example which shows that the appearance of $n_\Lambda$ is natural. Let $k$ denote the index for which $|\Lambda_k|=n_\Lambda$, and consider 
the index subsets
$$
\Lambda_{k,i}=\{\beta\in \Lambda_k:\;\beta_i=k\},
\qquad i=1,\ldots,d.
$$
Obviously, since $\cup_i \Lambda_{k,i}=\Lambda_k$, the largest of these subsets
has size at least $n_\Lambda/d$. Without loss of generality, we can assume that $\Lambda_{k,1}$ is this largest subset,
therefore the monotone index set 
$$
\Lambda':=\{\beta'\in \mathbb{N}^{d-1}:\; (k,\beta')\in \Lambda_{k,1}\}\subset\mathbb{N}^{d-1}
$$
satisfies $|\Lambda'|=|\Lambda_{k,1}|\ge n_\Lambda/d$. 

Consider now the function $\bar{s}_{\Lambda_{k,1}}\in S_\Lambda$ defined at the end of Section \ref{sec11}. According to (\ref{Lac1}), we have
$$
\|\bar{s}_{\Lambda_{k,1}}\|_{HB}^2 = 3^{-d} \sum_{\beta\in\Lambda_{k,1}} 2^{2|\beta|_\infty} =3^{-d} 2^{2k} |\Lambda'|.
$$
A lower bound for the $H_0^1$ norm of $\bar{s}_{\Lambda_{k,1}}$ is obtained as follows. Since
$$ 
\bar{s}_{\Lambda_{k,1}}(x) = \bar{s}_k(x_1) \bar{s}_{\Lambda'}(x'),
$$
by the construction of $\Lambda_{k,1}$, we get
$$
\left|\frac{\partial}{\partial x_1} \bar{s}_{\Lambda_{k,1}}(x)\right| = |\bar{s}'_k(x_1)| |\bar{s}_{\Lambda'}(x')|=2^k |\bar{s}_{\Lambda'}(x')|
$$
where we adopted the notation $x=(x_1,x')$ with $x'\in I^{d-1}$. Thus,
$$
\|\bar{s}_{\Lambda_{k,1}}\|_{H^1_0}^2 \ge \|\frac{\partial}{\partial x_1} \bar{s}_{\Lambda_{k,1}}\|_{L_2}^2
=2^{2k}\|\bar{s}_{\Lambda'}\|_{L_2}^2.
$$
Now we use (\ref{Lac2}) for $\bar{s}_{\Lambda'}$. This gives
$$
\|\bar{s}_{\Lambda_{k,1}}\|_{H^1_0}^2 \ge 4^{-d} 2^{2k}|\Lambda'|^2.
$$
Altogether we found that
\be\label{CLambda2}
C_\Lambda \ge \frac{\|\bar{s}_{\Lambda_{k,1}}\|_{H^1_0}^2}{\|\bar{s}_{\Lambda_{k,1}}\|_{HB}^2} \ge \frac{3^d}{4^d}|\Lambda'|
\ge \frac{3^d}{d4^d} n_\Lambda.
\ee
To summarize, according to (\ref{CLambda1}) and (\ref{CLambda2}) the best possible constant in (\ref{UE}) is proportional to $n_\Lambda$ up to constants only depending on $d$.

\section{Lower estimates: $c_{\Lambda}$}\label{sec3}
To estimate the best constant $c_\Lambda$ in (\ref{LE}), we proceed again as in \cite{GO1994} but start with a different decomposition
of an arbitrary $v_\Lambda\in S_\Lambda$. Namely, for $k=1,\ldots,k_\Lambda$, we define
the index subset $\Lambda_{k,0}\subset \Lambda_k$ by collecting into it all
maximal multi-indices in $\Lambda_k$, i.e., all $\beta\in \Lambda_k$ such that $\beta'\ge \beta$ and $\beta'\in \Lambda_k$ implies $\beta'=\beta$ (recall that $\Lambda_k$ and $k_\Lambda$ are given by (\ref{Lk}) ). 
Consider the $L_2$ orthogonal prewavelet representation
$$
v_\Lambda = \sum_{\beta\in\Lambda}  w_\beta =\sum_{k=1}^{k_\Lambda} \sum_{\beta'\in \Lambda_k} w_\beta',
$$
and partition, for each $k=1,\ldots,k_\Lambda$, the index set $\Lambda_k$ into subsets $\Lambda_{k,\beta}$ corresponding to the maximal indices $\beta\in\Lambda_{k,0}$ such that each $\beta'$ belongs to exactly one
$\Lambda_{k,\beta}$, and $\beta'\in \Lambda_{k,\beta}$ implies $\beta'\le \beta$ (this partitioning is always possible but not unique).  This gives a new decomposition
$$
v_{\Lambda}=\sum_{k=1}^{k_\Lambda}\sum_{\beta\in \Lambda_{k,0}} v_{\beta}, \qquad 
v_{\beta}:=\sum_{\beta'\in \Lambda_{k,\beta}} w_{\beta'}, \qquad \beta\in \Lambda_{k,0},
$$
into $v_\beta\in V_\beta$ for which 
$$
\|v_\beta\|_{L_2}^2=\sum_{\beta'\in \Lambda_{k,\beta}} \|w_{\beta'}\|_{L_2}^2
$$
due to the $L_2$ orthogonality of the $w_{\beta'}$.  Thus, according to (\ref{PWN}) we obtain the following lower bound for the $H_0^1$ norm of $v_\Lambda$:
\be\label{cLambda0}
\|v_\Lambda\|_{H_0^1}^2\ge c_{PW}\sum_{k=1}^{k_\Lambda} \sum_{\beta\in\Lambda_{k,0}} \sum_{\beta'\in \Lambda_{k,\beta}} \|w_{\beta'}\|_{L_2}^2=c_{PW} \sum_{k=1}^{k_\Lambda} \sum_{\beta\in\Lambda_{k,0}}  \|v_{\beta}\|_{L_2}^2.
\ee

We will next estimate the HB norms of the individual $v_\beta$ with $\beta\in \Lambda_{k,0}$
by their $L_2$ norms. 
For fixed but arbitrary $\beta\in \Lambda_{k,0}$, consider the HB decomposition 
$$
v_{\beta}= \sum_{\beta'\le \beta} s_{\beta'}.
$$
Each $s_{\beta'}$ is a telescoping sum of at most $2^d$ multi-linear spline interpolants $I_{\beta''}v_{\beta}\in V_{\beta''}$ of $v_{\beta}$
with respect to the nodal point set  $\Sigma_{\beta''}$ associated with tensor-product partition $T_{\beta''}$, where $\beta''\le \beta'$ and $|\beta'-\beta''|_\infty\le 1$.
Thus, the HB norm of $v_{\beta}$ is bounded by
\bea
\|v_{\beta}\|_{HB}^2 &=& \sum_{\beta'\le \beta} 2^{2|\beta'|_\infty} \|s_{\beta'}\|_{L_2}^2 \le
2^d \sum_{\beta'\le \beta} 2^{2|\beta'|_\infty} \sum_{\beta''\le \beta':\,|\beta'-\beta''|_\infty \le 1} \|I_{\beta''}v_{\beta}\|_{L_2}^2\\
&\le& C \sum_{\beta'\le \beta} 2^{2|\beta'|_\infty} \|I_{\beta'}v_{\beta}\|_{L_2}^2\le
C\sum_{\beta'\le \beta} 2^{2|\beta'|_\infty-|\beta'|_1} 
\sum_{P\in \Sigma_{\beta'}} |v_{\beta}(P)|^2.
\eea
Here and in the following we use the $L_2$-stability  
\be\label{L2S}
c\|v_{\beta}\|_{L_2}^2 \le 2^{-|\beta|_1} \sum_{P\in \Sigma_\beta} |v_\beta(P)|^2 \le C\|v_{\beta}\|_{L_2}^2,\qquad
v_\beta\in V_\beta,
\ee
of the nodal basis in the full grid spaces $V_\beta$ (in the above,  the lower $L_2$ stability  bound in (\ref{L2S})  was applied with $v_\beta$ replaced by  $I_{\beta'}v_\beta\in
V_{\beta'}$).
Since the nodal point sets $\Sigma_\beta$ form a monotone family with respect to the multi-index order, i.e., $\beta'\le \beta$ implies $\Sigma_{\beta'}\subset \Sigma_\beta$, 
we have
$$
\|v_{\beta}\|_{HB}^2\le C \sum_{P\in \Sigma_{\beta}} |v_{\beta}(P)|^2 
\sum_{\beta'\le \beta} 2^{2|\beta'|_\infty-|\beta'|_1}.
$$
A straightforward calculation shows that
$$
\sum_{\beta'\le \beta} 2^{2|\beta'|_\infty-|\beta'|_1}\le \sum_{i=1}^d \sum_{\beta'\le \beta:\,|\beta'|_\infty =\beta_i}
2^{2|\beta'|_\infty-|\beta'|_1}\le
d\sum_{k_1=1}^{k} 2^{k_1 }\sum_{1\le k_2,\ldots,k_d\le k_1} 2^{-k_2-\ldots-k_d} < d2^{k+1}.
$$
Thus,
\be\label{Help}
\|v_{\beta}\|_{HB}^2\le C 2^k \sum_{P\in \Sigma_\beta} |v_\beta(P)|^2 \le C2^{k+|\beta|_1} \|v_{\beta}\|_{L_2}^2,\qquad \beta\in\Lambda_{k,0},
\ee
for $k=1,\ldots,k_\Lambda$, where we have used the upper $L_2$ stability bound in (\ref{L2S}).

The combination of (\ref{cLambda0}) and  (\ref{Help}) leads to the bound
\bea
\|v_{\Lambda}\|_{HB}^2&\le& \|\sum_{k=1}^{k_\Lambda} \sum_{\beta\in\Lambda_{k,0}}
2^{(|\beta|_1-k)/2} 2^{(k-|\beta|_1)/2} v_\beta\|^2_{HB}\\
&\le& \left(\sum_{k=1}^{k_\Lambda} 2^{-k}\sum_{\beta\in\Lambda_{k,0}}
2^{|\beta|_1} \right)\left(\sum_{k=1}^{k_\Lambda}2^k \sum_{\beta\in\Lambda_{k,0}}
2^{-|\beta|_1} \|v_{\beta}\|_{HB}^2\right)  \\
&\le& C\tilde{n}_\Lambda \sum_{k=1}^{k_\Lambda}2^{2k} \sum_{\beta\in\Lambda_{k,0}} \|v_{\beta}\|_{L_2}^2    
\le C \tilde{n}_\Lambda \|v_{\Lambda}\|_{H^1_0}^2,
\eea
where
\be\label{NL1}
\tilde{n}_\Lambda:=\sum_{k=1}^{k_\Lambda} \sum_{\beta\in\Lambda_{k,0}}
2^{|\beta|_1-|\beta|_\infty}.
\ee
Thus, the best constant $c_\Lambda$ in (\ref{LE}) satisfies
\be\label{cLambda1}
c_\Lambda^{-1}\le C \tilde{n}_{\Lambda}, 
\ee
with $\tilde{n}_{\Lambda}$ defined in (\ref{NL1}).

The lower bound
\be\label{cLambda2}
c_\Lambda^{-1}\ge c\tilde{n}'_{\Lambda},\qquad \tilde{n}'_{\Lambda}:= \max_{\beta\in\Lambda} 2^{|\beta|_1-|\beta|_\infty}
\ee
is implied by considering the hat functions $\psi_\beta\in S_\Lambda$ with $\beta\in \Lambda$.
Indeed, (\ref{phiHB}) implies for each of them
$$
\|\psi_\beta\|_{HB}^2\ge c 2^{|\beta|_\infty},
$$
while  the same calculation that led to (\ref{phiN}) gives
$$
\|\psi_\beta\|_{H_0^1}^2\le C 2^{2|\beta|_\infty-|\beta|_1}.
$$
Thus,
$$
c_\Lambda^{-1}\ge \max_{\beta\in\Lambda} \frac{\|\psi_\beta\|_{HB}^2}{\|\psi_\beta\|_{H_0^1}^2} \ge c\max_{\beta\in\Lambda} 2^{|\beta|_1-|\beta|_\infty}
$$
yields (\ref{cLambda2}).

As we will see in the next section, for some important classes of generalized sparse grid spaces including $V_k$ and  $S_k$,
the two estimates (\ref{cLambda1}) and (\ref{cLambda2}) are of the same order but in general they do not match. Since $\Lambda_{k,0}\subset \Lambda_k\subset \Lambda$, we have
$$
1\le \frac{\tilde{n}_{\Lambda}}{\tilde{n}'_{\Lambda}} \le \sum_{k=1}^{k_\Lambda} \sum_{\beta\in\Lambda_{k,0}} 1 =
\sum_{k=1}^{k_\Lambda} |\Lambda_{k,0}| \le C\sum_{k=1}^{k_\Lambda} k^{d-2} \le Ck_\Lambda^{d-1}.
$$
Here we have used that $\Lambda_{k}$ is the union of $d$ sets $\Lambda_{k,i}$ each of which is essentially equivalent to 
a certain monotone set $\Lambda'\subset \mathbb{N}^{d-1}$ with $k_{\Lambda'}\le k$. Therefore, the cardinality of the set $\Lambda_{k,0}$
of maximal indices in $\Lambda_k$ cannot exceed $d$ times the maximal cardinality of the set of maximal indices of such $\Lambda'$ which
is bounded by $Ck^{d-2}$. 

Thus, in the worst case upper and lower bounds for $c_\Lambda^{-1}$ may be off by a factor $Ck_\Lambda^{d-1}$.
To see that such a gap is indeed possible, consider the index sets
$$
\Lambda = \{(2k,\beta'):\; \beta'\in \mathbb{N}^{d-1}, |\beta'|_1< k+d-1\}, \qquad k=1,2,\ldots,
$$
for which $k_\Lambda=2k$, and 
$$
\sum_{\beta\in \Lambda_{m,0}} 2^{|\beta|_1-|\beta|_\infty}=\sum_{(m,\beta')\in \Lambda_{m,0}} 2^{|\beta'|_1} = 2^k |\Lambda_{m,0}|\ge
ck^{d-2} 2^k,\qquad m=k+1,\ldots,2k.
$$
This gives $ \tilde{n}_{\Lambda}\ge ck^{d-1} 2^k$ and $\tilde{n}'_{\Lambda}=2^k$. We currently do not know if the gap can be reduced by constructing more sophisticated examples in order to improve the bound (\ref{cLambda2}).

\section{Summary}\label{sec4}
To summarize, we have established bounds for HB preconditioning of $H_0^1$-elliptic variational problems in generalized sparse grid 
spaces $S_\Lambda$ that are close to optimal for large classes of monotone $\Lambda$, in particular, for $V_k$,  $S_k$, and the energy-optimized
sparse grid spaces defined in \cite{GK2009}, see also \cite{GK2000,BG2004}.
\begin{theo} \label{theo1} Let $d>1$, and $\Lambda\subset \mathbb{N}^d$ be a monotone index set. The condition number $\kappa_{S_\Lambda,HB}$ of the tensor-product HB preconditioner 
of a discretization of a symmetric $H^1_0$-elliptic variational problem with the generalized sparse grid space $S_\Lambda$ 
satisfies the two-sided estimate
$$
c\, n_{\Lambda}\tilde{n}'_{\Lambda}\le \kappa_{S_\Lambda,HB}\le C\, n_{\Lambda}\tilde{n}_{\Lambda},
$$
where $n_{\Lambda}, \tilde{n}_{\Lambda}, \tilde{n}'_{\Lambda}$ are defined in (\ref{CLambda1}), (\ref{NL1}), (\ref{cLambda2}), respectively, and the constants $c,C$ depend solely on $d$.\\
For the standard sparse grid spaces $S_k$, the estimate turns into 
$$
\kappa_{S_k,HB} \approx k^{d-1} 2^{k(d-1)/d},\qquad k\to\infty.
$$
For the isotropic full grid spaces $V_k$, we have
$$
\kappa_{V_k,HB} \approx k^{d-1}2^{k(d-1)},\qquad k\to\infty.
$$
More generally, for the energy-optimized sparse-grid spaces $S_k^a:= S_{\Lambda_{k,a}}$, $-\infty < a < 1$, given by the index set
\be\label{Lka}
\Lambda_{k,a}:= \{\beta\in \mathbb{N}^d:\; |\beta|_1-a|\beta|_{\infty} \le (1-a)k + d-1\},
\ee
we have
\be\label{Ska}
\kappa_{S^a_k,HB} \approx k^{d-1} 2^{k(d-1)(1-a)/(d-a)},\qquad k\to\infty,
\ee
with constants that depend on $d$ (and may depend on $a$).
\end{theo}

Before proving the asymptotic condition number behavior for the families $V_k$, $S_k$, and $S_k^a$, let us mention that the result for $\kappa_{S_k,HB}$ improves our previous estimates for $d=2$ and $d=3$ stated in \cite{GO1994} by a factor $k$ and $k^3$, respectively. It is interesting to note that the condition number growth of the tensor-product HB preconditioner for $V_k$
is always worse compared to the isotropic HB preconditioner (papers by Yserentant \cite{Y1986,Y1993} for $d=2$, Ong \cite{O1997} for $d=3$, see 
\cite{O1994} for arbitrary $d>3$) by roughly an exponential factor of $2^k$. 

As to the energy-optimized sparse grid spaces $S_k^a$ which are defined in \cite[Section 4.1.2]{GK2009} as $V_J^T$ using different notation, it is obvious that $S_k=S_k^0$, $S_k^a$ becomes $V_k$ if $a\to -\infty$, and $S_k^a$ deteriorates for $a=1$ into a sum of essentially univariate spline spaces
$$
S_k^{1} = V_{k,1,\ldots,1} + V_{1,k,\ldots,1} +\ldots + V_{1,1,\ldots,k}.
$$

{\bf Proof} of Theorem \ref{theo1}. The condition number estimate for general $\Lambda$ is an immediate consequence of the results of Section \ref{sec2} and \ref{sec3}. The gap between upper and lower bounds is related to the quotient $\tilde{n}_{\Lambda} /\tilde{n}'_{\Lambda} $ which may grow as $k_\Lambda^{d-1}$ in the worst case.

The result for $V_k$ is obvious since for the associated $\Lambda$ we have $k_\Lambda=k$, and the sets $\Lambda_{r,0}$ consist of a single index $(r,\ldots,r)$,
$r=1,\ldots,k$. Therefore, 
$$
\tilde{n}'_{\Lambda} = 2^{k(d-1)},\qquad \tilde{n}_{\Lambda} = \sum_{r=1}^k 2^{r(d-1)} \le C2^{k(d-1)},
$$
while $n_\Lambda=|\Lambda_k| \approx k^{d-1}$. 

To prove (\ref{Ska}) for the generalized sparse grid spaces $S_k^a$ (which includes $S_k=S_k^0$ as a special case) we observe the following: As long as the integer $r$ satisfies the inequality $dr -ar \le (1-a)k + d-1$ or, equivalently, $r\le r_0:=\lfloor ((1-a)k + d-1)/(d-a)\rfloor$, we have
$V_r\subset S_k^a$ by the definition (\ref{Lka}) of the index set $\Lambda=\Lambda_k^a$. Thus, the sets $\Lambda_{r,0}$ of extremal points in $\Lambda_r$
consist of a single index $(r,r,\ldots,r)$ for $r=1,\ldots,r_0$ and 
$$
\sum_{r=1}^{r_0} \sum_{\beta\in\Lambda_{r,0}} 2^{|\beta|_1-|\beta|_\infty} = \sum_{r=1}^{r_0} 2^{(d-1)r}\le C2^{(d-1)r_0}.
$$
For the remaining $r=r_0+1,\ldots,k$ (note that $k_\Lambda=k$ since 
$$
(1-a)|\beta|_\infty + (d-1)\le |\beta|_1-a |\beta|_\infty\le k(1-a) +d-1, \qquad \beta\in \Lambda_k^a,
$$
and $(k,1,\ldots,1)\in \Lambda_k^a$ for any $a<1$), we split $\Lambda_r$ and thus $\Lambda_{r,0}$ into $d$ (not necessarily disjoint) sets
$\Lambda_r^i=\{\beta\in \Lambda_r:\;\beta_i=r\}$. Obviously, the maximal elements of $\Lambda_r^i$ the set of which is denoted $\Lambda^i_{r,0}$
belong to $\Lambda_{r,0}$, and, vice versa, each index in $\Lambda_{r,0}$ belongs to at least one $\Lambda^i_{r,0}$.
Since
$$
\sum_{\beta\in\Lambda_{r,0}} 2^{|\beta|_1-|\beta|_\infty}\le\sum_{i=1}^d \sum_{\beta\in\Lambda_{r,0}^i} 2^{|\beta|_1-|\beta|_\infty}=
d\sum_{\beta\in\Lambda_{r,0}^1} 2^{|\beta|_1-|\beta|_\infty}
$$
due to the invariance of $\Lambda_k^a$ with respect to index permutations, we estimate the sum for 
$$
\Lambda_{r,0}^1=\{(r,\beta'):\;\beta'\in \mathbb{N}^{d-1},\;|\beta'|_\infty\le r,\;|\beta'|_1=\lfloor (1-a)(k-r)+d-1 \rfloor\}.
$$
This description of  $\Lambda_{r,0}^1$ holds because any $\beta\in\Lambda_{r}^1$ is of the form $\beta=(r,\beta')$ with $|\beta'|_\infty\le r$ and satisfies
$$
|\beta|_1-a|\beta|_{\infty} = r(1-a) + |\beta'|_1\le (1-a)k + d-1,
$$
and therefore any maximal index $\beta$ for $\Lambda_{r}^1$ must satisfy equality $|\beta'|_1=\lfloor (1-a)(k-r)+d-1 \rfloor$.
This implies that
$$
\sum_{\beta\in\Lambda_{r,0}^1} 2^{|\beta|_1-|\beta|_\infty}=\sum_{\beta\in\Lambda_{r,0}^1} 2^{|\beta'|_1}=2^{\lfloor (1-a)(k-r)+d-1 \rfloor}
|\Lambda_{r,0}^1|.
$$
Note that $2^{\lfloor (1-a)(k-r)+d-1 \rfloor}\le C 2^{(1-a)(k-r_0)+d-1}2^{-(1-a)(r-r_0)}$, where the upper bound decays geometrically for any $a<1$. 
Thus, if we prove that
\be\label{Surf}
|\Lambda_{r,0}^1| \le C (r-r_0)^{d-2}, \qquad r=r_0+1,\ldots,k,
\ee
then 
$$
\sum_{r=r_0+1}^k \sum_{\beta\in\Lambda_{r,0}} 2^{|\beta|_1-|\beta|_\infty}\le C 2^{(1-a)(k-r_0)+d-1}\sum_{r=r_0+1}^k (r-r_0)^{d-2}2^{-(1-a)(r-r_0)}
\le C2^{(d-1)r_0},
$$
where we have used that $(d-1)r_0=(d-a)r_0 - (1-a)r_0\approx (1-a)k+d-1-(1-a)r_0=(1-a)(k-r_0)+d-1$ as $k\to\infty$. Thus, upper and lower bound for 
$c_{\Lambda_k^a}^{-1}$ are up to constant factors matching since
$$
\tilde{n}_{\Lambda_k^a} \le C2^{(d-1)r_0}\le C\tilde{n}'_{\Lambda_k^a}\le C\tilde{n}_{\Lambda_k^a},\qquad k=1,2,\ldots.
$$
(recall that $(r_0,r_0,\ldots,r_0)\in \Lambda_k^a$ and consequently $\tilde{n}'_{\Lambda_k^a}\ge 2^{(d-1)r_0}$). The result is
\be\label{cLka}
c_{\Lambda_k^a}^{-1} \approx 2^{(d-1)r_0} \le C2^{(d-1)(1-a)k/(d-a)},\qquad k=1,2,\ldots.
\ee
On the other hand, $C_{\Lambda_k^a}\approx n_{\Lambda_k^a}\approx k^{d-1}$ (we leave this to the reader). Together
this implies (\ref{Ska}).

It remains to show (\ref{Surf}). For $d=2$ this is obvious, since $\beta'\in \mathbb{N}$ and fixing $|\beta'|_1$ means to fix $\beta'$, i.e.,
in this case $|\Lambda_{r,0}^1|=1$ for all $r=r_0+1,\ldots,k$. For $d>2$, we set $\gamma_i=r-\beta'_i\ge 0$, $i=1,\ldots,d-1$,  and observe that
\bea
|\Lambda_{r,0}^1|&\le& |\{\gamma\in \mathbb{Z}_+^{d-1}:\; |\gamma|_1 = (d-1)r-\lfloor (1-a)(k-r)+d-1 \rfloor\}| \\
&\le& C((d-1)r-\lfloor (1-a)(k-r)+d-1 \rfloor)^{d-2}.
\eea
This follows from $|\beta'|_\infty\le r$ and $|\gamma|_1=(d-1)r-|\beta'|_1$ for all $\beta=(r,\beta')\in \Lambda_{r,0}^1$. Since
\bea
(d-1)r-\lfloor (1-a)(k-r)+d-1 \rfloor &=& (d-a)r - ((1-a)k+d-1) +\epsilon \\
&=& (d-a)(r-r_0) +\epsilon -\epsilon'(d-a),
\eea
where $0\le \epsilon,\epsilon'<1$, we have 
$$
|(d-1)r-\lfloor (1-a)(k-r)+d-1 \rfloor|\le C(r-r_0), \qquad r=r_0+1,\ldots,k,
$$ 
with a constant $C$ depending on $d$ and $a$. This gives 
(\ref{Surf}) and concludes the proof of Theorem \ref{theo1}.\hfill $\Box$

\medskip
We mention that in \cite[Section 3.2]{BG2004} a modification of $\Lambda_k^{1/5}$ is used to define a so-called energy-based sparse grid space, in order to optimize error bounds in the $H_0^1$ norm, and claim that the method used for (\ref{Ska}) covers this example as well. 

Finding the correct order of the HB condition numbers $\kappa_{V_\beta,HB}$ for arbitrary anisotropic full grid spaces
$V_\beta$ as $\beta \to \infty$ is currently an open problem. By inspecting the subsets $\Lambda_{k,0}$ for the index set $\Lambda$ associated with $V_\beta$ reveals that here the gap between $\tilde{n}'_{\Lambda}$
and $\tilde{n}_{\Lambda}$ may be as large as  $k_\Lambda=|\beta|_\infty$ (but not larger), independently of $d$. 
Indeed, without loss of generality, set $\beta=(k_\Lambda,\beta')$ with $\beta'\in \mathbb{N}^{d-1}$ satisfying $|\beta'|_\infty\le k_\Lambda$. Then $\Lambda_{k,0}=\{(k,\min(\beta',(k,\ldots,k))\}$, $k=1,\ldots,k_\Lambda$, where the minimum of the  two index vectors is taken componentwise. Thus, 
$$
\tilde{n}'_{\Lambda}=\max_{k=1,\ldots,k_\Lambda} 2^{|\min(\beta',(k,\ldots,k))|_1} = 2^{|\beta'|_1},
$$
while
$$
\tilde{n}_{\Lambda}=\sum_{k=1}^{k_\Lambda} 2^{|\min(\beta',(k,\ldots,k))|_1} \le k_\Lambda 2^{|\beta'|_1}.
$$
In the last estimate, equality is attained for $\beta'=(1,\ldots,1)$. A direct inspection of this extreme case of 
$V_\beta$ with $\beta=(k_\Lambda,1,\ldots,1)$ shows that $c_\Lambda^{-1}\approx 1$ independently of 
$k_\Lambda$ and $d$. In other words, in this case the lower bound (\ref{cLambda2}) gives the  correct behavior which indicates that improvements in the proof of the upper bound  (\ref{cLambda1}) for
$c_\Lambda^{-1}$ should be possible.

\section*{Acknowledgement} 
This paper was written during the author's stay at the Institute for Numerical Simulation (INS) sponsored by the
Hausdorff Center for Mathematics of the University of Bonn funded by the Deutsche Forschungsgemeinschaft. He is grateful for this support and the stimulating atmosphere at the INS.

\end{document}